\documentclass[12pt]{article}%
\usepackage{amsfonts,amsmath}
\usepackage{amsmath}
\usepackage{amsfonts}
\usepackage{amssymb}
\usepackage{graphicx}%
\providecommand{\U}[1]{\protect\rule{.1in}{.1in}}
\begin{document}

\title{{A closed formula for the topological entropy of multimodal maps based on
min-max symbols}}
\author{Jos\'{e} M. Amig\'{o}\thanks{Corresponding author. E-mail: jm.amigo@umh.es.
Phone: +34966658911.}, and \'{A}ngel Gim\'{e}nez\\Centro de Investigaci\'{o}n Operativa, Universidad Miguel Hern\'{a}ndez.\\Avda. de la Universidad s/n. 03202 Elche (Alicante). Spain.}
\maketitle
\date{}

\vspace*{0.5cm}

\noindent\textbf{Abstract}\quad Topological entropy is a measure of complex
dynamics. In this regard, multimodal maps play an important role when it comes
to study low-dimensional chaotic dynamics or explain some features of higher
dimensional complex dynamics with conceptually simple models. In the first
part of this paper an analytical formula for the topological entropy of twice
differentiable multimodal maps is derived, and some basic properties are
studied. This expression involves the so-called min-max symbols, which are
closely related to the kneading symbols. Furthermore, its proof leads to a
numerical algorithm that simplifies a previous one also based on min-max
symbols. In the second part of the paper this new algorithm is used to compute
the topological entropy of different modal maps. Moreover, it compares
favorably to the previous algorithm when computing the topological entropy of
the bi- and tri-modal maps considered in the numerical simulations.

\section{Introduction}

The kneading sequences of a multimodal map $f$ are symbolic sequences that
locate the iterates of its critical values up to the precision set by the
partition defined by its critical points \cite{Milnor,Melo}. Since the $n$th
iterate of a critical point of $f$ is a critical value of $f$, it makes sense
to attach to the symbols of each kneading sequence of $f$ a label informing
about their minimum/maximum (or \textquotedblleft critical\textquotedblright)
character. The result is a \textit{min-max sequence}, one per critical point,
consisting of min-max symbols. These symbols and sequences were introduced in
\cite{Dias1,Dilao} for unimodal maps, and in \cite{Amigo2012} for multimodal
maps. Thus, min-max sequences generalize kneading sequences in that they
additionally provide geometric information about the extrema structure of
$f^{n}$ at the critical points for all $n\geq1$. That this generalization is a
good idea, can be justified in several ways, the most direct one being that
the computational cost of a min-max symbol is virtually the same as of a
kneading symbol for any sufficiently smooth multimodal map. Indeed, the extra
piece of information contained in a min-max symbol can be automatically
retrieved from a look-up table once the min-max symbol of the previous iterate
has been calculated.

Another justification is that min-max sequences allow to construct recursive
algorithms to compute the topological entropy \cite{Adler,Walters} of
multimodal maps. To this end we assumed in \cite{Dilao2012,Amigo2012} that $f$
is twice differentiable, although numerical simulations with continuous,
piecewise linear maps of constant slopes $\pm\left\vert s\right\vert $ (and
hence, with topological entropy $\log\left\vert s\right\vert $ \cite{Alseda})
support the hypothesis that our results hold true under weaker conditions.

In this paper, which is an outgrowth of \cite{Amigo2012}, we derive an
analytical formula for the topological entropy of $f$, $h(f)$, that is
formally similar to other well-known expressions like
\cite{Misiu,Alseda,Balm}
\begin{align}
h(f)  &  =\lim_{n\rightarrow\infty}\frac{1}{n}\log\ell_{n}\label{1.1}\\
&  =\lim_{n\rightarrow\infty}\frac{1}{n}\log\left\vert \left\{  x\in
I:f^{n}(x)=x\right\}  \right\vert \label{1.2}\\
&  =\lim_{n\rightarrow\infty}\frac{1}{n}\log^{+}\text{Var}(f^{n})\label{1.3}\\
&  =\lim_{n\rightarrow\infty}\frac{1}{n}\log^{+}\text{length}(f^{n}),
\label{1.4}%
\end{align}
where (i) $\ell_{n}$ is shorthand for the \textit{lap number} of $f^{n}$
(i.e., the number of maximal monotonicity segments of $f^{n}$), (ii)
$\left\vert \cdot\right\vert $ denotes cardinality (i.e., $\left\vert \left\{
x\in I:f^{n}(x)=x\right\}  \right\vert $ is the number of periodic points of
period $n$), (iii) Var$(f^{n})$ stands for the variation of $f^{n}$, and (iv)
length$(f^{n})$ means the length of the graph of $f^{n}$. The new expression
follows from (\ref{1.1}) via some geometrical properties for boundary-anchored
maps involving min-max symbols. Moreover, its derivation leads to several
numerical algorithms to compute $h(f)$. We will only discuss the most simple
one, which abridges the algorithm of \cite{Amigo2012}. Benchmarking of the
simplified algorithm with respect to the original one shows that the former
sometimes outperforms the latter.

The interest of closed (or analytical) formulas is manifold including, as just
mentioned, hypothetical improvements in the speed and precision of already
existing algorithms. But, most importantly, analytical formulas usually
provide insights into a problem, open alternative ways to prove new properties
or attack old problems and, in any case, add techniques to the conceptual and
instrumental toolkit of the field.\textbf{ }Thus, as compared to the general
definition of topological entropy \cite{Adler,Walters}, the expressions
(\ref{1.1})-(\ref{1.4}) are conceptually simpler, besides providing a variety
of numerical techniques to compute $h(f)$; see \cite{Dilao2012,Amigo2012} for
general algorithms based on the formula (\ref{1.1}), and
\cite{Balm,Baldwin,Block,Block2,Collet,Froyland,Gora,Steinberger} for other
mathematical schemes with various degrees of generality.

This paper is organized as follows. In order to make the paper self-contained,
we review in Sect.~\ref{section2} all the basic concepts, especially the
concept of min-max sequences, needed for the present sequel. In
Sect.~\ref{section3} we introduce some auxiliary results. In particular, we
provide a formal proof of the known fact that the topological entropy does not
depend on the boundary conditions. Sect. \ref{section4} contains the main
result of the paper, namely, a new analytical formula for the topological
entropy of multimodal maps (Theorem 1). As way of illustration, this formula
is applied in Sect.~\ref{section5} to a few special cases of multimodal maps
whose critical values comply with certain confinement conditions. In
Sect.~\ref{section6} we derive an interesting relation between the value of
$h(f)$ and the divergence rate of a logarithmic expression that appears in the
analytical formula proved in Theorem 1. A simple algorithm prompted by the
proof of Theorem 1 is explained in Sect. \ref{section7}. This algorithm is put
to test in Sect.~\ref{section8}, where the topological entropy of uni-, bi-,
and trimodal maps taken from \cite{Dilao2012,Amigo2012} are computed. Its
performance is also compared with the algorithm of \cite{Amigo2012} in those
three cases.

\section{Min-max sequences}

\label{section2}

We use the same notation as in \cite{Amigo2012} throughout.

Let $I$ be a compact interval $[a,b]\subset\mathbb{R}$ and $f:I\rightarrow I$
a piecewise monotone continuous map. Such a map is called $l$-modal if $f$ has
precisely $l$ \textit{turning points} (i.e., points in $(a,b)$ where $f$ has a
local extremum). Assume that $f$ has local extrema at $c_{1}<...<c_{l}$ and
that $f$ is strictly monotone in each of the $l+1$ intervals%
\[
I_{1}=[a,c_{1}),I_{2}=(c_{1},c_{2}),...,I_{l}=(c_{l-1},c_{l}),I_{l+1}%
=(c_{l},b]\text{.}%
\]
Also as in \cite{Amigo2012}, we assume henceforth that $f(c_{1}$) is a
maximum. These maps are said to have positive shape. This implies that
$f(c_{2k+1})$, $0\leq k\leq\left\lfloor \frac{l-1}{2}\right\rfloor $, are
maxima, whereas $f(c_{2k})$, $1\leq k\leq\left\lfloor \frac{l}{2}\right\rfloor
$, are minima. Furthermore, $f$ is strictly increasing on the the intervals
$I_{2k+1}$, $0\leq k\leq\left\lfloor \frac{l}{2}\right\rfloor $, and strictly
decreasing on the intervals $I_{2k}$, $1\leq k\leq\left\lfloor \frac{l+1}%
{2}\right\rfloor $.

The \textit{itinerary} of $x\in I$ under $f$ is a symbolic sequence
\[
\mathbf{i}(x)=(i_{0}(x),i_{1}(x),...,i_{n}(x),...)\in\{I_{1},c_{1}%
,I_{2},...,c_{l},I_{l+1}\}^{\mathbb{N}_{0}}%
\]
($\mathbb{N}_{0}\equiv\{0\}\cup\mathbb{N}$), defined as follows:%
\[
i_{n}(x)=\left\{
\begin{array}
[c]{cl}%
I_{j} & \text{if }f^{n}(x)\in I_{j}\text{ }(1\leq j\leq l+1),\\
c_{k} & \text{if }f^{n}(x)=c_{k}\text{ }(1\leq k\leq l).
\end{array}
\right.
\]
The itineraries of the critical values,%

\[
\gamma^{i}=(\gamma_{1}^{i},...,\gamma_{n}^{i},...)=\mathbf{i}(f(c_{i}%
)),\;1\leq i\leq l,
\]
are called the \textit{kneading sequences \cite{Melo}} of $f$.

We turn to the min-max sequences \cite{Dias1,Dilao,Amigo2012,Dilao2012} of $f$.

\medskip

\noindent\textbf{Definition 1.} The \textit{min-max sequences} of an $l$-modal
map $f$,
\[
\omega^{i}=(\omega_{1}^{i},\omega_{2}^{i},...,\omega_{n}^{i},...),\;1\leq
i\leq l,
\]
are defined as follows:%
\[
\omega_{n}^{i}=\left\{
\begin{array}
[c]{ll}%
m^{\mathbf{\gamma}_{n}^{i}} & \text{if }f^{n}(c_{i})\text{ is a minimum},\\
M^{\mathbf{\gamma}_{n}^{i}} & \text{if }f^{n}(c_{i})\text{ is a maximum.}%
\end{array}
\right.
\]
where $\gamma_{n}^{i}$ are kneading symbols.

Thus, the \textit{min-max} \textit{symbols} $\omega_{n}^{i}$ have an
exponential-like notation, where the `base' belongs to the alphabet $\{m,M\}$,
and the `exponent' is a kneading symbol. Therefore, the extra information of a
min-max symbol $\omega^{i}$ as compared to a kneading symbol $\gamma^{i}$ lies
in the base.

As in \cite{Dilao2012,Amigo2012} we consider the class $\mathcal{F}_{l}(I)$ of
$l$-modal maps $f:I\rightarrow I$ such that

(i) $f\in C^{2}(I)$, and

(ii) $f^{\prime}(x)\neq0$ for $x\in I_{k}$, $1\leq k\leq l+1$.

\noindent That is, a multimodal map $f$ of the interval $I$ belongs to
$\mathcal{F}_{l}(I)$ if (i) it is twice differentiable, and (ii) it is
strictly monotone except at the turning points. When the interval $I$ is clear
from the context or unimportant for the argument, we write just $\mathcal{F}%
_{l}$.

The next lemma proves our claim in the Introduction that, from the point of
view of the computational cost, min-max sequences and kneading sequences are
virtually equivalent, at least if $f$ is twice differentiable.

\medskip

\noindent\textbf{Lemma 1.} If $f\in\mathcal{F}_{l}$ has positive shape, then
the following `transition rules' hold:%
\begin{equation}%
\begin{tabular}
[c]{|c|c|c|}\hline
$\omega_{n}^{i}$ & $\rightarrow$ & $\omega_{n+1}^{i}$\\\hline\hline
$m^{c_{\text{even}}},M^{c_{\text{even}}}$ & $\rightarrow$ & $m^{\gamma
_{n+1}^{i}}$\\\hline
$m^{c_{\text{odd}}},M^{c_{\text{odd}}}$ & $\rightarrow$ & $M^{\gamma_{n+1}%
^{i}}$\\\hline
$m^{I_{\text{odd}}},M^{I_{\text{even}}}$ & $\rightarrow$ & $m^{\gamma
_{n+1}^{i}}$\\\hline
$m^{I_{\text{even}}},M^{I_{\text{odd}}}$ & $\rightarrow$ & $M^{\gamma
_{n+1}^{i}}$\\\hline
\end{tabular}
\ \ \ \ \ \label{transition}%
\end{equation}
where \textquotedblleft even\textquotedblright\ and \textquotedblleft
odd\textquotedblright\ stand for even and odd subindices, respectively, of the
critical points $c_{1},...,c_{l}$, and of the intervals $I_{1},...,I_{l}$.

\medskip

See \cite[Lemma 2.2]{Amigo2012}. Therefore, the kneading symbols of the
$f\in\mathcal{F}_{l}$, together with its \textit{initial min-max symbols},
i.e.%
\begin{equation}
\omega_{1}^{i}=\left\{
\begin{array}
[c]{ll}%
M^{\gamma_{1}^{i}} & \text{if }i=1,3,...,2\left\lfloor \frac{l+1}%
{2}\right\rfloor -1,\\
m^{\gamma_{1}^{i}} & \text{if }i=2,4,...,2\left\lfloor \frac{l}{2}%
\right\rfloor ,
\end{array}
\right.  \;\label{initial}%
\end{equation}
and the transition rules (\ref{transition}) allow to compute the min-max
sequences of $f\in\mathcal{F}_{l}$ in a recursive way.

In \cite{Dilao2012,Amigo2012} we used `signatures' instead of kneading symbols
in the exponents of $\omega_{n}^{i}$. The \textit{signature} of a point
$x\in\lbrack a,b]$ is a vector with $l$ entries, the $i$th entry being $+1$,
$0$, or $-1$ according to whether $x>c_{i}$, $x=c_{i}$, or $x<c_{i}$,
respectively. It is clear that the signature of $f^{n}(c_{i})$ does the same
as $\gamma_{n}^{i}$ when it comes to locate $f^{n}(c_{i})$ in the partition%
\begin{equation}
I_{1}\cup\{c_{1}\}\cup I_{2}\cup\{c_{2}\}\cup...,\cup\{c_{l}\}\cup I_{l+1}
\label{partition}%
\end{equation}
of the interval $I=[a,b]$, but in a `computer-friendly' way. For the purposes
of this paper though, the computational advantages of signatures will be not needed.

A final ingredient (proper of min-max sequences) is the following. Let the
$i$\textit{th critical line}, $1\leq i\leq l$, be the line $y=c_{i}$ in the
Cartesian product $I\times I=\{(x,y):x,y\in I\}$. Min-max symbols divide into
\textit{bad} and \textit{good symbols} with respect to $i$th critical line.
Geometrically, good symbols correspond to local maxima strictly above the line
$y=c_{i}$, or to local minima strictly below the line $y=c_{i}$. All other
min-max symbols are bad by definition with respect to the $i$th critical line.
We use the notation%
\[
\mathcal{B}^{i}=\{M^{I_{1}},M^{c_{1}},...,M^{I_{i}},M^{c_{i}},m^{c_{i}%
},m^{I_{i+1}},...,m^{c_{l}},m^{I_{l+1}}\}
\]
for the set of bad symbols of $f\in\mathcal{F}_{l}$ with respect to the $i$th
critical line. There are $2(l+1)$ bad symbols and $2(l-1)$ good symbols with
respect to a given critical line.

Bad symbols appear in all results of \cite{Dilao2012,Amigo2012} concerning the
computation of the topological entropy of $f\in\mathcal{F}_{l}$ via min-max
symbols. In this sense we may say that bad symbols are the hallmark of this approach.

\section{Auxiliary results}

\label{section3}

In general, Latin indices refer to the critical points and range between $1$
and $l$, while Greek indices refer to the number of iterations, hence they
take on arbitrary, nonnegative integer values.

Let $s_{\nu}^{i}$, $1\leq i\leq l$, stand for the \textit{number of interior
simple zeros of} $f^{\nu}(x)-c_{i}$, $\nu\geq0$, i.e., solutions of
$x-c_{i}=0$ ($\nu=0$), or (ii) solutions of $f^{\nu}(x)=c_{i}$, $x\in(a,b)$,
with $f^{\mu}(x)\neq c_{i}$ for $0\leq\mu\leq\nu-1$, and $f^{\nu\prime
}(x)\not =0$ ($\nu\geq1$). Geometrically $s_{\nu}^{i}$ is the number of
transversal intersections in the Cartesian plane $(x,y)$ of the curve
$y=f^{\nu}(x)$ and the straight line $y=c_{i}$, over the interval $(a,b)$.
Note that $s_{0}^{i}=1$ for all $i$.

To streamline the notation of the forthcoming math, set%
\begin{equation}
s_{\nu}=\sum_{i=1}^{l}s_{\nu}^{i} \label{snu}%
\end{equation}
for $\nu\geq0$. In particular,
\begin{equation}
s_{0}=\sum_{i=1}^{l}s_{0}^{i}=\sum_{i=1}^{l}1=l. \label{seed33}%
\end{equation}
According to \cite[Eqn. (31)]{Amigo2012}, the lap number of $f^{n}$, $\ell
_{n}$, satisfies
\begin{equation}
\ell_{n}=\sum\limits_{\nu=0}^{n-1}s_{\nu}+1. \label{ln-main}%
\end{equation}
In particular, $\ell_{1}=l+1$.

Furthermore, define%
\begin{equation}
K_{\nu}^{i}=\{(k,\kappa),1\leq k\leq l,1\leq\kappa\leq\nu:\omega_{\kappa}%
^{k}\in\mathcal{B}^{i}\}, \label{Kni}%
\end{equation}
($\nu\geq1$, $1\leq i\leq l$), that is, $K_{\nu}^{i}$ collects the upper and
lower indices $(k,\kappa)$ of the \textit{bad} symbols with respect to the
$i$th critical line in all the initial segments%
\[
\omega_{1}^{1},\omega_{2}^{1},...,\omega_{\nu}^{1};\;\;\omega_{1}^{2}%
,\omega_{2}^{2},...,\omega_{\nu}^{2};\;\;...;\;\;\omega_{1}^{l},\omega_{2}%
^{l},...,\omega_{\nu}^{l};
\]
of the min-max sequences of $f$. We note for further reference that $K_{\nu
-1}^{i}\subset K_{\nu}^{i}$, the set-theoretical difference being%
\begin{equation}
K_{\nu}^{i}\backslash K_{\nu-1}^{i}=\{(k,\nu),1\leq k\leq l:\omega_{\nu}%
^{k}\in\mathcal{B}^{i}\}. \label{Kni2}%
\end{equation}

Finally, set%
\begin{equation}
S_{\nu}^{i}=2\sum_{(k,\kappa)\in K_{\nu}^{i}}s_{\nu-\kappa}^{k},
\label{notation0}%
\end{equation}
where $S_{\nu}^{i}=0$ if $K_{\nu}^{i}=\emptyset$, and analogously to
(\ref{snu}),
\begin{equation}
S_{\nu}=\sum_{i=1}^{l}S_{\nu}^{i}. \label{notation}%
\end{equation}

The algorithm to compute the topological entropy of $f\in\mathcal{F}_{l}$ in
\cite{Amigo2012} rests on the relation \cite[Eqn. (32)]{Amigo2012}%
\begin{equation}
s_{\nu}^{i}=1+\sum_{\mu=0}^{\nu-1}s_{\mu}-S_{\nu}^{i}-\alpha_{\nu}^{i}%
-\beta_{\nu}^{i}, \label{sinu0}%
\end{equation}
where $\alpha_{\nu}^{i}$,$\beta_{\nu}^{i}$ are binary variables that depend on
$f^{\nu}(a)$, $f^{\nu}(b)$, and the $i$th critical line in the way specified
in \cite[Eqn. (27)]{Amigo2012}. Let us point out for further reference that
all $\alpha_{\nu}^{i}$'s and $\beta_{\nu}^{i}$'s vanish if $f$ is
\textit{boundary-anchored}, i.e., $f\{a,b\}\subset\{a,b\}$. Since we are
considering $l$-modal maps with a positive shape, this condition boils down in
our case to%
\[
f(a)=f(b)=a
\]
if $l$ is odd, or%
\[
f(a)=a,\;f(b)=b
\]
if $l$ is even.

We are going to show that, as long as the computation of $h(f)$ is concerned,
we may assume without loss of generality that $f$ is boundary-anchored. Its
proof proceeds by extending $f$ to a selfmap $F$ of a greater interval in such
a way that $F$ is boundary-anchored, and $h(f)=h(F)$.

To prove this property, we need the following general facts. Let
$g:X\rightarrow X$ be a continuous map of a compact Hausdorff space $X$ into
itself. A point $x\in X$ is nonwandering with respect to the map $g$ if for
any neighborhood $U$ of $x$ there an $n\geq1$ (possibly depending on $x$) such
that $f^{n}(U)\cap U\neq\emptyset$. Fixed and periodic points are examples of
nonwandering points. The closed set of all nonwandering points of $g$ is
called its \textit{nonwandering set} and denoted by $\Omega(g)$. According to
\cite[Lemma 4.1.5]{Alseda},%
\begin{equation}
h(g)=h(\left.  g\right\vert _{\Omega(g)}). \label{fact1}%
\end{equation}
Furthermore, if%
\[
X=\bigcup\limits_{i=1}^{k}Y_{i}%
\]
and all $Y_{i}$ are closed and $g$-invariant (i.e., $g(Y_{i})\subset Y_{i})$,
then \cite[Lemma 4.1.10]{Alseda},%
\begin{equation}
h(g)=\max_{1\leq i\leq k}h(\left.  g\right\vert _{Y_{i}}). \label{fact2}%
\end{equation}

\noindent\textbf{Lemma 2.} Let $f\in\mathcal{F}_{l}(I)$. Then there exists
$F\in\mathcal{F}_{l}(J)$, where $J\supset I$, such that $h(F)=h(f)$ and $F$ is boundary-anchored.

\medskip

\noindent\textit{Proof}. Set $I=[a,b]$, and $J=[a^{\prime},b^{\prime}]$ with
$a^{\prime}\leq a<b\leq b^{\prime}$. If $f(a)=a$, choose $a^{\prime}=a$; if
$f(b)=a$ ($l$ odd) or $f(b)=b$ ($l$ even), choose $b^{\prime}=b$. For
definiteness, we suppose the most general situation, namely, $a^{\prime}<a$
\ and $b<b^{\prime}$. Let $F:J\rightarrow J$ be such that (i) $F$ is strictly
increasing and twice differentiable on $[a^{\prime},a]$, (ii) $\left.
F\right\vert _{[a,b]}=f$, and (iii) $F$ is strictly decreasing ($l$ odd) or
strictly increasing ($l$ even), and twice differentiable on $[b,b^{\prime}]$.
Moreover the extension of $f$ to $F$ can be made in such a way that $F$ is
twice differentiable at the points $a$ and $b$, hence $F\in\mathcal{F}_{l}(J)$
by construction. As a result, $F$ has the same critical points and values as
$f$, and it is boundary-anchored.

Furthermore it is easy to check that $\Omega(F)=\Omega(f)\cup C$, where $C$ is
a closed and $F$-invariant set that only contains fixed points. Thus,
$h(\left.  F\right\vert _{C})=0$ and, according to (\ref{fact1}) and
(\ref{fact2}),%
\[
h(F)=h(\left.  F\right\vert _{\Omega(F)})=\max\{h(\left.  F\right\vert
_{\Omega(f)}),h(\left.  F\right\vert _{C})\}=h(\left.  F\right\vert
_{\Omega(f)})=h(\left.  f\right\vert _{\Omega(f)})=h(f).\;\square
\]

\medskip

The formulation and proof of Lemma 2 was tailored to maps $f\in\mathcal{F}%
_{l}(I)$ of positive shape. It is plain though that the statement of Lemma 2
holds also true if $f$ is just a continuous selfmap of a closed interval $I$.
In this case, $F$ may be taken piecewise linear on $[a^{\prime},a]\cup\lbrack
b,b^{\prime}]$.

\section{A closed formula for the topological entropy of unimodal maps}

\label{section4}

According to Lemma 2, given $f\in\mathcal{F}_{l}$ we may assume without
restriction that it is boundary-anchored when calculating its topological
entropy. This being the case, set $\alpha_{\nu}^{i}=\beta_{\nu}^{i}=0$ in
(\ref{sinu0}) for all $\nu\geq1$ and $i=1,...,l$, i.e.,
\begin{equation}
s_{\nu}^{i}=1+\sum_{\mu=0}^{\nu-1}s_{\mu}-S_{\nu}^{i}, \label{sinu}%
\end{equation}
and sum (\ref{sinu}) over $i$ from $1$ to $l$ to obtain the relation%

\begin{equation}
s_{\nu}=l\left(  \sum\limits_{\mu=0}^{\nu-1}s_{\mu}+1\right)  -S_{\nu}
\label{account}%
\end{equation}
between $s_{0}=l,s_{1},...,s_{\nu}$ and $S_{\nu}$, for all $\nu\geq1$. By
(\ref{ln-main}) this equation can we rewritten as $s_{\nu}=l\ell_{\nu}-S_{\nu
}$, hence%
\begin{equation}
\ell_{\nu}=\frac{1}{l}(S_{\nu}+s_{\nu}). \label{account3}%
\end{equation}

\medskip

\noindent\textbf{Lemma 3.} Let $f\in\mathcal{F}_{l}$ be boundary-anchored.
Then%
\begin{equation}
s_{\nu}=l(l+1)^{\nu}-l\sum_{\delta=1}^{\nu-1}(l+1)^{\nu-\delta-1}S_{\delta
}-S_{\nu} \label{induction}%
\end{equation}
for $\nu\geq1$, where the summation over $\delta$ is missing for $\nu=1$.

\medskip

\noindent\textit{Proof}. The proof is by induction. The case $\nu=1$ holds
trivially on account of (\ref{account}) and $s_{0}=l$.

Consider next the case $\nu+1$. By (\ref{account}) with $\nu+1$ instead of
$\nu$,
\begin{align}
s_{\nu+1}  &  =l+l\sum\limits_{\mu=0}^{\nu}s_{\mu}-S_{\nu+1}=l+ls_{0}%
+l\sum\limits_{\mu=1}^{\nu}s_{\mu}-S_{\nu+1}\nonumber\\
&  =l(1+l)+l\sum\limits_{\mu=1}^{\nu}\left(  l(l+1)^{\mu}-l\sum_{\delta
=1}^{\mu-1}(l+1)^{\mu-\delta-1}S_{\delta}-S_{\mu}\right)  -S_{\nu
+1}\label{inductionh1}\\
&  =l(l+1)^{\nu+1}-l\sum\limits_{\mu=1}^{\nu}\left(  l\sum_{\delta=1}^{\mu
-1}(l+1)^{\mu-\delta-1}S_{\delta}+S_{\mu}\right)  -S_{\nu+1} \label{3terms1}%
\end{align}
The induction hypothesis (\ref{induction}) was applied in line
(\ref{inductionh1}). The middle term in (\ref{3terms1}) can be simplified as
follows.%
\begin{align*}
\sum\limits_{\mu=1}^{\nu}\left(  l\sum_{\delta=1}^{\mu-1}(l+1)^{\mu-\delta
-1}S_{\delta}+S_{\mu}\right)   &  =l\sum\limits_{\mu=1}^{\nu}\sum_{\delta
=1}^{\mu-1}(l+1)^{\mu-\delta-1}S_{\delta}+\sum\limits_{\mu=1}^{\nu}S_{\mu}\\
&  =l\sum_{\delta=1}^{\nu-1}\sum_{k=0}^{\nu-\delta-1}(l+1)^{k}S_{\delta}%
+\sum\limits_{\mu=1}^{\nu}S_{\mu}\\
&  =\sum_{\delta=1}^{\nu-1}\left(  (l+1)^{\nu-\delta}-1\right)  S_{\delta
}+\sum\limits_{\mu=1}^{\nu}S_{\mu}\\
&  =\sum_{\delta=1}^{\nu-1}(l+1)^{\nu-\delta}S_{\delta}+S_{\nu}\\
&  =\sum_{\delta=1}^{\nu}(l+1)^{\nu-\delta}S_{\delta}.
\end{align*}
Replacement of this into (\ref{3terms1}) yields%
\begin{equation}
s_{\nu+1}=l(l+1)^{\nu+1}-l\sum_{\delta=1}^{\nu}(l+1)^{\nu-\delta}S_{\delta
}-S_{\nu+1}. \label{induction2}%
\end{equation}
Comparison of (\ref{induction2}) with (\ref{induction}) completes the
induction step. $\;\square$

\medskip

\noindent\textbf{Theorem 1.} Let $f\in\mathcal{F}_{l}$. Then%
\begin{equation}
h(f)=\log(l+1)-\lim_{\nu\rightarrow\infty}\frac{1}{\nu}\log\frac{1}%
{1-\sum_{\delta=1}^{\nu-1}\frac{S_{\delta}}{(l+1)^{\delta+1}}},
\label{Theorem}%
\end{equation}
with $S_{\delta}$ as in (\ref{notation}).

\medskip

\noindent\textit{Proof}. Without restriction, suppose that $f$ is boundary
anchored. From (\ref{account3}) and (\ref{induction}),
\begin{equation}
\ell_{\nu}=\frac{S_{\nu}+s_{\nu}}{l}=(l+1)^{\nu}-\sum_{\delta=1}^{\nu
-1}(l+1)^{\nu-\delta-1}S_{\delta}=(l+1)^{\nu}\left(  1-\sum_{\delta=1}^{\nu
-1}\frac{S_{\delta}}{(l+1)^{\delta+1}}\right)  . \label{Theoremb}%
\end{equation}
Use now (\ref{1.1}) to derive%
\[
h(f)=\lim_{\nu\rightarrow\infty}\frac{1}{\nu}\log\ell_{\nu}=\log
(l+1)+\lim_{\nu\rightarrow\infty}\frac{1}{\nu}\log\left(  1-\sum_{\delta
=1}^{\nu-1}\frac{S_{\delta}}{(l+1)^{\delta+1}}\right)  .\;\square
\]

Note that (\ref{Theoremb}) along with (\ref{sinu})-(\ref{induction}) hold true
only for \textit{boundary-anchored} $f\in\mathcal{F}_{l}$. In other words:
while the lap numbers $\ell_{\nu}$ depend in general both on the min-max
sequences of $f$ (through the sets $K_{\nu}^{i}$ in (\ref{notation0})) and the
itineraries of the boundary points of $I$ (through the $\alpha_{\nu}^{i}$'s
and $\beta_{\nu}^{i}$'s in (\ref{sinu0})), the limit $\lim_{\nu\rightarrow
\infty}\frac{1}{\nu}\log\ell_{\nu}=h(f)$ only depends on the min-max sequences.

According to (\ref{Theorem}), $h(f)\leq\log(l+1)$, a well-known result for
multimodal maps. Moreover, (\ref{Theorem}) expresses the difference
$\log(l+1)-h(f)$ for $f\in\mathcal{F}_{l}$ with the help of $(S_{\delta
})_{\delta\geq1}$. The computation of $h(f)$, based on (\ref{account3}), will
be addressed in Sect. \ref{section7} and \ref{section8}.

\section{Special cases}

\label{section5}

As way of illustration of the results of Sect.~\ref{section3}, let us consider
a few special cases characterized by fulfilling what we call `confinement
conditions' for the critical values.

\medskip

\textbf{(C1)} If $l$ is \textit{odd} and $f(c_{i})<c_{1}$ for $i=1,...,l$ (so
that the graphs of $y=f^{n}(x)$ lie below the first critical line $y=c_{1}$),
then (see (\ref{transition})) $\omega^{\text{odd}}=(M^{I_{1}})^{\infty}$
contains only bad symbols with respect to all critical lines, and
$\omega^{\text{even}}=(m^{I_{1}})^{\infty}$ contains only good symbols with
respect to all critical lines. Thus,
\[
K_{\nu}^{i}=\bigcup\limits_{k=0}^{(l-1)/2}\{(2k+1,1),(2k+1,2),...,(2k+1,\nu
)\}.
\]
Furthermore, $s_{0}^{i}=1$, and $s_{n}^{i}=0$ for $n\geq1$, $1\leq i\leq l$,
in this case. It follows%
\[
S_{\nu}=2\sum_{i=1}^{l}\sum_{(k,\kappa)\in K_{\nu}^{i}}s_{\nu-\kappa}%
^{k}=2\sum_{i=1}^{l}(s_{0}^{1}+s_{0}^{3}+...+s_{0}^{l})=2l\frac{(l+1)}%
{2}=l(l+1),
\]
hence
\begin{align*}
h(f)  &  =\log(l+1)+\lim_{\nu\rightarrow\infty}\frac{1}{\nu}\log\left(
1-\sum_{\delta=1}^{\nu-1}\frac{l}{(l+1)^{\delta}}\right) \\
&  =\log(l+1)+\lim_{\nu\rightarrow\infty}\frac{1}{\nu}\log\frac{1}%
{(l+1)^{\nu-1}}=0.
\end{align*}

\textbf{(C2) }If $l$ is \textit{even} and $f(c_{i})<c_{1}$ for all
$i=1,...,l$, then the min-max sequences are the same as in case (C1),
\[
K_{\nu}^{i}=\bigcup\limits_{k=0}^{l/2-1}\{(2k+1,1),(2k+1,2),...,(2k+1,\nu)\},
\]
but this time, due to the branch of $y=f^{n}(x)$ connecting $(c_{l}%
,f^{n}(c_{l}))$ with $(b,b)$, we have $s_{n}^{i}=1$ for $n\geq0$. It follows%
\[
S_{\nu}=2\sum_{i=1}^{l}\sum_{(k,\kappa)\in K_{\nu}^{i}}s_{\nu-\kappa}%
^{k}=2\sum_{i=1}^{l}\sum_{k=0}^{l/2-1}(s_{\nu-1}^{2k+1}+s_{\nu-2}%
^{2k+1}+...+s_{0}^{2k+1})=2l\tfrac{l}{2}\nu=l^{2}\nu,
\]
hence%
\begin{align*}
h(f)  &  =\log(l+1)+\lim_{\nu\rightarrow\infty}\frac{1}{\nu}\log\left(
1-l^{2}\sum_{\delta=1}^{\nu-1}\frac{\delta}{(l+1)^{\delta+1}}\right) \\
&  =\log(l+1)+\lim_{\nu\rightarrow\infty}\frac{1}{\nu}\log\frac{l\nu
+1}{(l+1)^{\nu}}=0.
\end{align*}

\textbf{(C3)} If $l$ \textit{even} and $f(c_{i})>c_{l}$ for all $i=1,...,l$,
then $\omega^{\text{odd}}=(M^{I_{l+1}})^{\infty}$ and $\omega^{\text{even}%
}=(m^{I_{l+1}})^{\infty}$. In this case,%
\[
K_{\nu}^{i}=\bigcup\limits_{k=1}^{l/2}\{(2k,1),(2k,2),...,(2k,\nu)\}.
\]
Analogously to case (C2) we also have $s_{n}^{i}=1$ for $n\geq0$, but now
owing to the branch of $y=f^{n}(x)$ connecting $(a,a)$ with $(c_{1}%
,f^{n}(c_{1}))$. Thus%
\[
S_{\nu}=2\sum_{i=1}^{l}\sum_{(k,\kappa)\in K_{\nu}^{i}}s_{\nu-\kappa}%
^{k}=2\sum_{i=1}^{l}\sum_{k=1}^{l/2}(s_{\nu-1}^{2k}+s_{\nu-2}^{2k}%
+...+s_{0}^{2k})=2l\tfrac{l}{2}\nu=l^{2}\nu,
\]
as in (C2), hence%
\[
h(f)=0\text{.}%
\]

\textbf{(C4)} Finally if $f(c_{\text{odd}})=b$, and $f(c_{\text{even}})=a$,
then (see (\ref{transition}))%
\[
\omega^{\text{odd}}=\left\{
\begin{array}
[c]{ll}%
M^{I_{l+1}}(m^{I_{1}})^{\infty} & \text{if }l\text{ is odd,}\\
(M^{I_{l+1}})^{\infty} & \text{if }l\text{ is even,}%
\end{array}
\right.
\]
and%
\[
\omega^{\text{even}}=(m^{I_{1}})^{\infty}\text{.}%
\]
Thus, in either case both $\omega^{\text{odd}}$ and $\omega^{\text{even}}$
contain only good symbols with respect to all critical lines. It follows that
$K_{\nu}^{i}=\emptyset$, hence $S_{\nu}=0$ for all $\nu\geq1$. We conclude%
\begin{equation}
h(f)=\log(l+1)\text{,} \label{hmax}%
\end{equation}
which is the maximum value $h(f)$ can achieve for $f\in\mathcal{F}_{l}$.

\section{Convergence}

\label{section6}

Let us study next the convergence of (\ref{Theorem}). From (\ref{Theoremb}) it
follows%
\begin{equation}
1-\sum_{\delta=1}^{\nu-1}\frac{S_{\delta}}{(l+1)^{\delta+1}}=\frac{\ell_{\nu}%
}{(l+1)^{\nu}}\geq0 \label{Theoremc}%
\end{equation}
for all $\nu\geq2$. Hence%
\begin{equation}
h(f)=\log(l+1)-\lim_{\nu\rightarrow\infty}\frac{1}{\nu}\log\frac{1}%
{1-\Sigma_{\nu-1}}, \label{hA}%
\end{equation}
where (see (\ref{Theoremc}))%
\begin{equation}
0\leq\Sigma_{\nu}:=\sum_{\delta=1}^{\nu}\frac{S_{\delta}}{(l+1)^{\delta+1}%
}\leq1 \label{Anu}%
\end{equation}
for $\nu\geq1$. By definition (\ref{Anu}), $\Sigma_{1},..,\Sigma_{\nu}%
,\Sigma_{\nu+1},...$ is a non-negative, non-decreasing sequence of real
numbers bounded by 1. Therefore, it converges and
\[
\Sigma_{\infty}:=\sum_{\delta=1}^{\infty}\frac{S_{\delta}}{(l+1)^{\delta+1}%
}=\lim_{\nu\rightarrow\infty}\Sigma_{\nu}\in\lbrack0,1].
\]

As way of example consider the special cases of Sect.~\ref{section5}. In case
(C1)%
\[
\Sigma_{\infty}=\sum_{\delta=1}^{\infty}\frac{l}{(l+1)^{\delta}}=1,
\]
in cases (C2) and (C3)%
\[
\Sigma_{\infty}=\sum_{\delta=1}^{\infty}\frac{l^{2}\delta}{(l+1)^{\delta+1}%
}=1,
\]
while $\Sigma_{\infty}=0$ in case (C4).

Remember that $A(n)=o(B(n))$ means that $\lim_{n\rightarrow\infty}%
A(n)/B(n)=0$, and $A(n)\sim B(n)$ means $\lim_{n\rightarrow\infty}A(n)/B(n)=1$.

\medskip

\noindent\textbf{Theorem 2. }Let $f\in\mathcal{F}_{l}$. Then%
\begin{equation}
h(f)\in\left\{
\begin{array}
[c]{ll}%
\{\log(l+1)\} & \text{if }\log(1-\Sigma_{n})=o(n),\\
\lbrack0,\log(l+1)) & \text{if }\log(1-\Sigma_{n})\sim-Cn,
\end{array}
\right.  \label{inclusion}%
\end{equation}
where $0<C\leq\log(l+1)$.

\medskip

\noindent\textit{Proof}. According to (\ref{hA}),
\[
h(f)=\log(l+1)+\lim_{n\rightarrow\infty}\frac{1}{n}\log(1-\Sigma_{n-1}).
\]
Therefore, $h(f)=\log(l+1)$ if $\lim_{n\rightarrow\infty}\frac{1}{n}%
\log(1-\Sigma_{n-1})=0$, i.e., $\log(1-\Sigma_{n})=o(n)$. Otherwise, $0\leq
h(f)<\log(l+1)$ if $\lim_{n\rightarrow\infty}\frac{1}{n}\log(1-\Sigma
_{n-1})=-C$ with $0<C\leq\log(l+1)$, i.e., $\log(1-\Sigma_{n})\sim
-Cn$.$\;\square$

\medskip

\medskip

\begin{figure}[h]
\centering
\includegraphics[width=10cm]{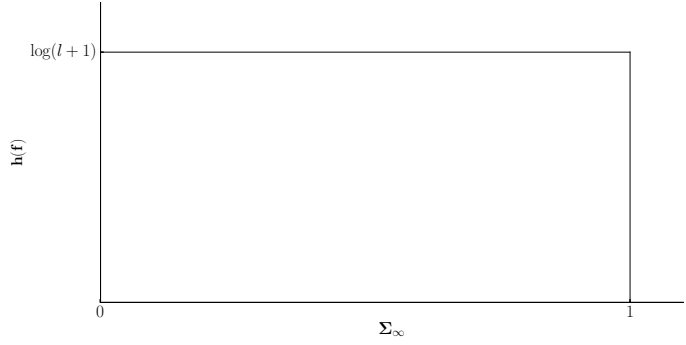} \caption{Graph of the inclusion
(\ref{inclusion3}).}%
\label{fig1}%
\end{figure}

\medskip

\noindent\textbf{Theorem 3. }Let $f\in\mathcal{F}_{l}$. Then%
\begin{equation}
h(f)\in\left\{
\begin{array}
[c]{ll}%
\{\log(l+1)\} & \text{if }0\leq\Sigma_{\infty}<1,\\
\lbrack0,\log(l+1)] & \text{if }\Sigma_{\infty}=1.
\end{array}
\right.  \label{inclusion3}%
\end{equation}

\medskip

\noindent\textit{Proof}. Taking the limit $n\rightarrow\infty$ in
(\ref{inclusion}), one finds that the correspondence $\Sigma_{\infty}\mapsto
h(f)$ defines the following inclusion:%
\begin{equation}
h(f)\in\left\{
\begin{array}
[c]{ll}%
\{\log(l+1)\} & \text{if }0\leq\Sigma_{\infty}<1,\text{ or }\Sigma_{\infty
}=1\text{ with }\log(1-\Sigma_{n})=o(n),\\
\lbrack0,\log(l+1)) & \text{if }\Sigma_{\infty}=1\text{ with }\log
(1-\Sigma_{n})\sim-Cn.
\end{array}
\right.  \label{inclusion2}%
\end{equation}
This proves (\ref{inclusion3}).$\;\square$

\medskip

Note that in the second case of (\ref{inclusion2}), $h(f)=\log(l+1)-C$ (see
the proof of Theorem 2) and $\Sigma_{n}\sim1-e^{-Cn}$, hence%
\[
\sum_{\delta=1}^{n}\frac{S_{\delta}}{(l+1)^{\delta+1}}\sim1-e^{-\left(
\log(l+1)-h(f)\right)  n}.
\]
In other terms, $\Sigma_{n}\nearrow1$ exponentially fast when $h(f)<\log
(l+1)$, the difference $1-\Sigma_{n}$ decreasing as $e^{-\left(
\log(l+1)-h(f)\right)  n}$.

Figure 1 depicts the inclusion (\ref{inclusion3}): If $h(f)=\log(l+1)$ then
$0\leq\Sigma_{\infty}\leq1$; if $\Sigma_{\infty}=1$ then $0\leq h(f)\leq
\log(l+1)$.

\section{\textbf{A simplified algorithm for the topological entropy}}

\label{section7}

When it comes to calculate numerically $h(f)=\lim_{\nu\rightarrow\infty}%
\frac{1}{\nu}\log\ell_{\nu}$ via Eqn. (\ref{Theoremb}), used in the proof of
Theorem 1, the intermediate expression%
\begin{equation}
h(f)=\lim_{\nu\rightarrow\infty}\frac{1}{\nu}\log\frac{S_{\nu}+s_{\nu}}{l}
\label{efficient}%
\end{equation}
is more efficient and numerically stable than the final expression
\begin{equation}
h(f)=\log(l+1)+\lim_{\nu\rightarrow\infty}\frac{1}{\nu}\log\left(
1-\sum_{\delta=1}^{\nu-1}\frac{S_{\delta}}{(l+1)^{\delta+1}}\right)  .
\label{Theorem2}%
\end{equation}
The computation of $S_{\nu}^{i}$, $1\leq i\leq l$, requires $s_{0}^{i}%
=1,s_{1}^{i},...,s_{\nu-1}^{i}$, see (\ref{notation0}), while the computation
of $s_{\nu}^{i}$, $1\leq i\leq l$, requires $s_{0}^{i},s_{1}^{i},...,s_{\nu
-1}^{i}$, and $S_{\nu}^{i}$, see (\ref{sinu}).

We summarize the algorithm resulting from (\ref{efficient}) in the following
scheme (\textquotedblleft$A\longrightarrow B$\textquotedblright\ stands for
\textquotedblleft$B$ is computed by means of $A$\textquotedblright).

\begin{description}
\item[(A1)] \textbf{Parameters: }$l\geq1$ (number of critical points),
$\varepsilon>0$ (dynamic halt criterion), and $n_{\max}\geq2$ (maximum number
of loops).

\item[(A2)] \textbf{Initialization:} $s_{0}^{i}=1$, and $K_{1}^{i}=\{k,1\leq
k\leq l:\omega_{1}^{k}\in\mathcal{B}^{i}\}$ ($1\leq i\leq l$).

\item[(A3)] \textbf{First iteration: }For $1\leq i\leq l$,\textbf{ }%
\[%
\begin{array}
[c]{rcll}%
s_{0}^{i},K_{1}^{i} & \longrightarrow & S_{1}^{i},S_{1} & \text{(use
(\ref{notation0}), (\ref{notation}))}\\
s_{0}^{i},S_{1}^{i} & \longrightarrow & s_{1}^{i},s_{1} & \text{(use
(\ref{sinu}), (\ref{account}))}%
\end{array}
\]

\item[(A4)] \textbf{Computation loop}. For $1\leq i\leq l$ and $\nu\geq2$ keep
calculating $K_{\nu}^{i}$, $S_{\nu}^{i}$, and $s_{\nu}^{i}$ according to the
recursions
\begin{equation}%
\begin{array}
[c]{rcll}%
K_{\nu-1}^{i} & \longrightarrow & K_{\nu}^{i} & \text{(use (\ref{Kni2}),
(\ref{transition}))}\\
s_{0}^{i},s_{1}^{i},...,s_{\nu-1}^{i},K_{\nu}^{i} & \longrightarrow & S_{\nu
}^{i},S_{\nu} & \text{(use (\ref{notation0}), (\ref{notation}))}\\
s_{0}^{i},s_{1}^{i},...,s_{\nu-1}^{i},S_{\nu}^{i} & \longrightarrow & s_{\nu
}^{i},s_{\nu} & \text{(use (\ref{sinu}), (\ref{account}))}%
\end{array}
\label{loop}%
\end{equation}
until (i)
\begin{equation}
\left\vert \frac{1}{\nu}\log\frac{S_{\nu}+s_{\nu}}{l}-\frac{1}{\nu-1}\log
\frac{S_{\nu-1}+s_{\nu-1}}{l}\right\vert \leq\varepsilon, \label{halt}%
\end{equation}
or, else, (ii) $\nu=n_{\max}+1$.

\item[(A5)] \textbf{Output.} In case (i) output%
\begin{equation}
h(f)=\frac{1}{\nu}\log\frac{S_{\nu}+s_{\nu}}{l}. \label{approx}%
\end{equation}
In case (ii) output \textquotedblleft Algorithm failed\textquotedblright.
\end{description}

The algorithm (A1)-(A5) simplifies the original algorithm \cite{Amigo2012},
which is based on the exact value of the lap number $\ell_{\nu}$. This entails
that the new algorithm needs more loops to output $h(f)$ with the same
parameter $\varepsilon$ in the halt criterion (\ref{halt}), although this does
not necessarily mean that the overall execution time will be longer since now
less computations are required. In fact, we will find both situations in the
numerical simulations below.

Two final remarks:

\begin{description}
\item[R1.] The parameter $\varepsilon$ does not bound the error $\left\vert
h(f)-\frac{1}{\nu}\log\frac{S_{\nu}+s_{\nu}}{l}\right\vert $ but the
difference between two consecutive estimations, see (\ref{halt}). The number
of exact decimal positions of $h(f)$ can be found out by taking different
$\varepsilon$'s , as we will see in the next section. Equivalently, one can
control how successive decimal positions of $\frac{1}{\nu}\log\frac{S_{\nu
}+s_{\nu}}{l}$ stabilize with growing $\nu$. Moreover, the smaller $h(f)$, the
smaller $\varepsilon$ has to be chosen to achieve a given approximation precision.

\item[R2.] According to \cite[Thm. 4.2.4]{Alseda}, $\frac{1}{\nu}\log\ell
_{\nu}\geq h(f)$ for any $\nu$. We may expect therefore that the numerical
approximations (\ref{approx}) converge from above to the true value of the
topological entropy with ever more iterations, in spite of the relation
$\ell_{\nu}=\frac{1}{l}(S_{\nu}+s_{\nu})$ holding in general for
boundary-anchored maps only.
\end{description}

\section{Numerical simulations}

\label{section8}

In this section we calculate the topological entropy of families of uni-, bi-,
and trimodal maps taken from \cite{Dilao2012} and \cite{Amigo2012}; none of
these maps is boundary-anchored. The purpose of our choice is to compare the
performance of the simplified algorithm with the original one. To this end, a
code for arbitrary $l$ was written with PYTHON, and run on an Intel(R)
Core(TM)2 Duo CPU. All logarithms were taken to base $e$, i.e., the values of
the topological entropy in this section are given in nats per iteration. The
numerical results will be given with six decimal positions for brevity.

\subsection{Simulation with 1-modal maps}

\label{section8.1}

Let $\alpha>0$, $-1<\beta\leq0$, and $f_{\alpha,\beta}:[-(1+\beta
),(1+\beta)]\rightarrow\lbrack-(1+\beta),(1+\beta)]$ be defined as \cite[Eqn.
(29)]{Dilao2012}
\[
f_{\alpha,\beta}(x)=e^{-\alpha^{2}x^{2}}+\beta.
\]
These maps have the peculiarity of showing direct and reverse period-doubling
bifurcations when the parameters are monotonically changed \cite[Fig.
3(a)]{Dilao2012}.

Fig. 2 shows the plot of $h(f_{2.8,\beta})$ vs $\beta$ calculated with the
algorithm of Sect. 7. Here $\varepsilon=10^{-4}$ and the parameter $\beta$ was
increased in steps of $\Delta\beta=0.001$ from $\beta=-0.999$ to $\beta=0$.
Upon comparing Fig. 2 with Fig. 3(b) of \cite{Dilao2012}, we see that both
plots coincide visually, except for the two vanishing entropy tails. We
conclude that $\varepsilon=10^{-4}$ is not small enough to obtain reliable
estimations of the topological entropy for vanishing values of $h(f_{2.8,\beta
})$. This fact can also be ascertained numerically by taking different values
of $\varepsilon$, as we do in the table below.

\begin{figure}[h]
\centering
\includegraphics[width=10cm]
	{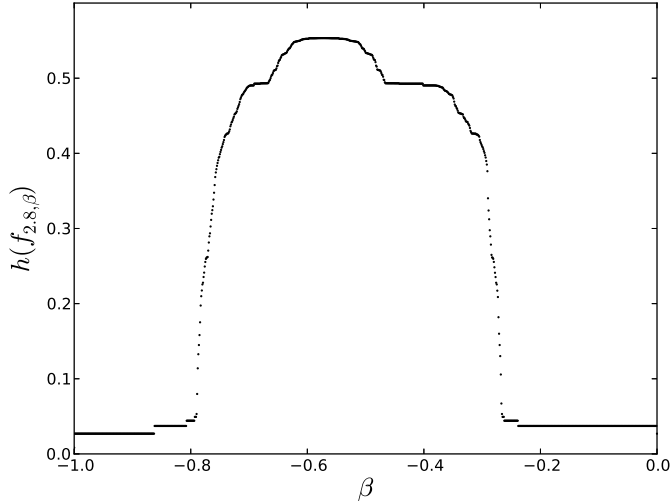}\caption{Plot of $h(f_{2.8,\beta})$ vs $\beta$, $-1<\beta\leq0$
($\varepsilon=10^{-4},\Delta\beta=0.001$).}%
\label{fig1a}%
\end{figure}

In order to compare the convergence speed and execution time of the original
(\cite{Dilao2012,Amigo2012}) and the simplified algorithm, we have computed
$h(f_{2.8,-0.5})$ with both algorithms for different $\varepsilon$'s. The
number of loops $n$ needed to achieve the halt condition $\varepsilon=10^{-d}%
$, $4\leq d\leq7$, and the execution time $t$ (in seconds) are listed in Table
1. The columns $h_{orig}$, $n_{orig}$, and $t_{orig}$ were obtained with the
original algorithm, while the columns $h_{simp}$, $n_{simp}$, and $t_{simp}$
were obtained with the simplified one. For all choices of $\varepsilon$,
$t_{orig}<t_{simp}$. Furthermore, we conclude from Table 1 that
$h(f)=0.524...$ using the original algorithm, while $h(f)=0.52...$ using the
simplified one and the same set of $\varepsilon$'s.

\begin{table}[h]
\centering
\begin{tabular}
[c]{|c|c|c|c|c|c|c|}\hline
& $h_{orig}$ & $n_{orig}$ & $t_{orig}$ & $h_{simp}$ & $n_{simp}$ & $t_{simp}%
$\\\hline
$\varepsilon=10^{-4}$ & 0.531968 & 81 & 0.024452 & 0.534106 & 101 &
0.023519\\\hline
$\varepsilon=10^{-5}$ & 0.526645 & 253 & 0.200861 & 0.527305 & 318 &
0.217069\\\hline
$\varepsilon=10^{-6}$ & 0.524935 & 797 & 1.87643 & 0.525142 & 1004 &
2.126501\\\hline
$\varepsilon=10^{-7}$ & 0.524391 & 2519 & 18.404195 & 0.524456 & 3174 &
21.399207\\\hline
\end{tabular}
\caption{Comparison of performances when computing $h(f_{2.8,-0.5})$.}%
\label{table1}%
\end{table}

Fig. 3 depicts the values of $h(f_{\alpha,\beta})$ for $2\leq\alpha\leq3$,
$-1<\beta\leq0$, $\varepsilon=10^{-4}$, and $\Delta\alpha,\Delta\beta=0.01$.

\begin{figure}[h]
\centering
\includegraphics[width=10cm]
	{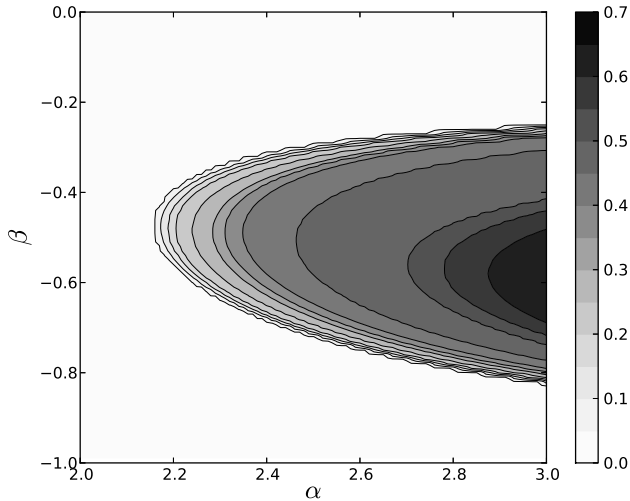}\caption{Level sets of $h(f_{\alpha,\beta})$ vs $\alpha,\beta$,
$2\leq\alpha\leq3$, and $-1<\beta\leq0$ ($\varepsilon=10^{-4},\Delta
\alpha=\Delta\beta=0.01$).}%
\label{fig1b}%
\end{figure}

\subsection{Simulation with 2-modal maps}

\label{section8.2}

Let $0\leq v_{2}<v_{1}\leq1$ and $f_{v_{1},v_{2}}:[0,1]\rightarrow\lbrack0,1]$
be defined as \cite[Sect. 8.1]{Amigo2012}
\[
f_{v_{1},v_{2}}(x)=(v_{1}-v_{2})(16x^{3}-24x^{2}+9x)+v_{2},
\]
These maps have convenient properties for numerical simulations as they share
the same fixed critical points,%
\[
c_{1}=1/4,\;c_{2}=3/4,
\]
the critical values are precisely the parameters,%
\[
f_{v_{1},v_{2}}(1/4)=v_{1},\;f_{v_{1},v_{2}}(3/4)=v_{2},
\]
and the values of $f$ at the endpoints are explicitly given by the parameters
as follows:%
\[
f_{v_{1},v_{2}}(0)=v_{2},\;f_{v_{1},v_{2}}(1)=v_{1}.
\]

Fig. 4 shows the plot of $h(f_{1,v_{2}})$ vs $v_{2}$, $0\leq v_{2}<1$,
computed with the new algorithm, $\varepsilon=10^{-4}$, and $\Delta
v_{2}=0.001$. Again, this plot coincides visually with the same plot computed
with the old algorithm \cite[Fig. 4]{Amigo2012} except for the vanishing
entropy tail, which indicates that $\varepsilon=10^{-4}$ is too large a value
for obtaining accurate estimates in that parametric region. \begin{figure}[h]
\centering
\includegraphics[width=10cm]
	{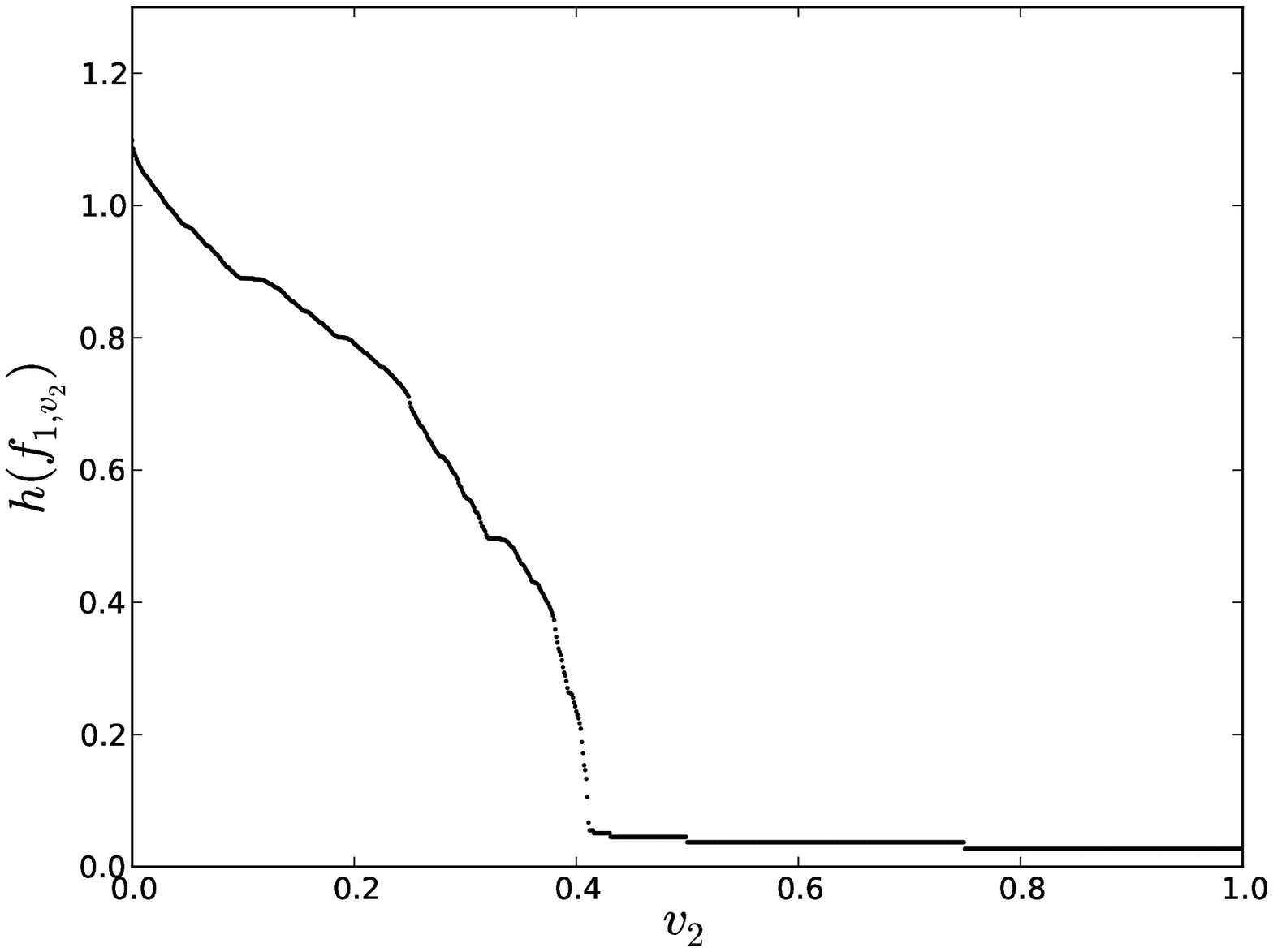}\caption{Plot of $h(f_{1,v_{2}})$ vs $v_{2}$, $0\leq v_{2}\leq1$
($\varepsilon=10^{-4},\Delta v_{2}=0.001$).}%
\end{figure}

Table 2 displays the performance of the new algorithm as compared to the old
one when computing $h(f_{0.9,0.1})$. This time $t_{orig}>t_{simp}$ for
$\varepsilon=10^{-d}$, $4\leq d\leq7$ (as in Table 1). Furthermore, we obtain
two exact decimal positions of the topological entropy, $h(f_{0.9,0.1}%
)=0.41...$, with both algorithms and the $\varepsilon$'s considered.

\begin{table}[h]
\centering
\begin{tabular}
[c]{|c|c|c|c|c|c|c|}\hline
& $h_{orig}$ & $n_{orig}$ & $t_{orig}$ & $h_{simp}$ & $n_{simp}$ & $t_{simp}%
$\\\hline
$\varepsilon=10^{-4}$ & 0.432246 & 162 & 0.228102 & 0.434175 & 182 &
0.202511\\\hline
$\varepsilon=10^{-5}$ & 0.421287 & 511 & 2.1436 & 0.42191 & 573 &
1.973616\\\hline
$\varepsilon=10^{-6}$ & 0.417812 & 1613 & 20.918574 & 0.418008 & 1810 &
19.738665\\\hline
$\varepsilon=10^{-7}$ & 0.41671 & 5099 & 207.190258 & 0.416772 & 5721 &
182.951993\\\hline
\end{tabular}
\caption{Comparison of performances when computing $h(f_{0.9,0.1})$.}%
\label{table2}%
\end{table}

Fig. 5 depicts the values of $h(f_{v_{1},v_{2}})$ for $0\leq v_{2}\leq
v_{1}-0.5$, $\varepsilon=10^{-4}$, and $\Delta v_{1},\Delta v_{2}=0.01$.

\begin{figure}[h]
\centering
\includegraphics[width=10cm]
	{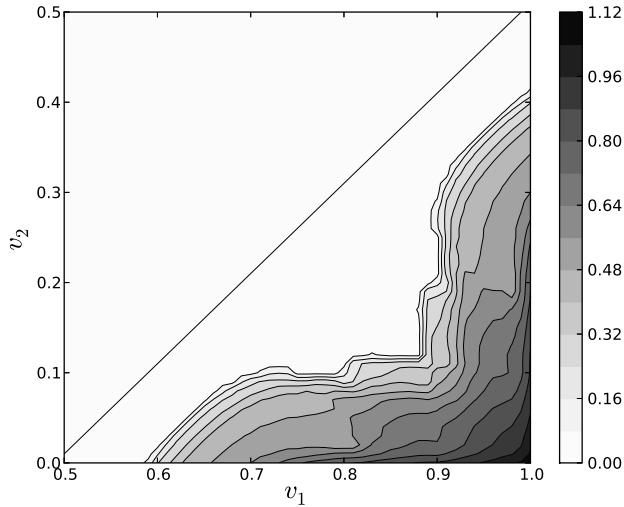}\caption{Level sets of $h(f_{v_{1},v_{2}})$ vs $v_{1},v_{2}$,
$0\leq v_{2}\leq v_{1}-0.5$ ($\varepsilon=10^{-4},\Delta v_{1}=\Delta
v_{2}=0.01$).}%
\label{fig2b}%
\end{figure}

\subsection{Simulation with 3-modal maps}

\label{section8.3}

Consider next the $3$-modal maps $f_{v_{2},v_{3}}:[0,1]\rightarrow\lbrack0,1]$
defined by the quartic polynomials \cite[Sect. 8.2]{Amigo2012}%
\begin{align*}
f_{v_{2},v_{3}}(x)  &  =\frac{4\left[  \left(  2\sqrt{2}-1\right)
v_{2}-2v_{3}\right]  x}{2(2\sqrt{2}+1)v_{3}-7v_{2}}\left[  4\left(
1+2\sqrt{2}\right)  (x-1)(1-2x)^{2}v_{3}\right. \\
&  +\left.  \left(  -56x^{3}+20\left(  4+\sqrt{2}\right)  x^{2}-\left(
37+18\sqrt{2}\right)  x+3\sqrt{2}+5\right)  v_{2}\right]  ,
\end{align*}
where $0\leq v_{2}<v_{3}\leq1$. The critical points of $f_{v_{2},v_{3}}$ are
\[
c_{1}=\frac{-\sqrt{2}\text{$v_{2}$}-4\text{$v_{2}$}+12\sqrt{2}\text{$v_{3}$%
}-8\text{$v_{3}$}}{8\left(  -7\text{$v_{2}$}+4\sqrt{2}\text{$v_{3}$%
}+2\text{$v_{3}$}\right)  },\quad c_{2}=1/2,\quad c_{3}=\frac{1}{4}(2+\sqrt
{2}).
\]
Moreover this family verifies $f_{v_{2},v_{3}}(0)=0$, $f_{v_{2},v_{3}}%
(c_{2})=v_{2}$, $f(c_{3})=v_{3}$, and%

\[
f_{v_{2},v_{3}}(1)=\frac{4\left(  5\sqrt{2}-8\right)  \text{$v_{2}$}\left(
\left(  2\sqrt{2}-1\right)  \text{$v_{2}$}-2\text{$v_{3}$}\right)
}{-7\text{$v_{2}$}+4\sqrt{2}\text{$v_{3}$}+2\text{$v_{3}$}}.
\]

Fig. 6 shows the plot of $h(f_{v_{2},1})$ vs $v_{2}$, $0\leq v_{2}<1$,
computed with the new algorithm, $\varepsilon=10^{-4}$, and $\Delta
v_{2}=0.001$. Once more, this plot coincides visually with the same plot
computed with the old algorithm \cite[Fig. 7 (left)]{Amigo2012} except for the
vanishing entropy tail, which again indicates that $\varepsilon=10^{-4}$ is
too large a value for obtaining accurate estimates in that parametric region.

\begin{figure}[h]
\centering
\includegraphics[width=10cm]
	{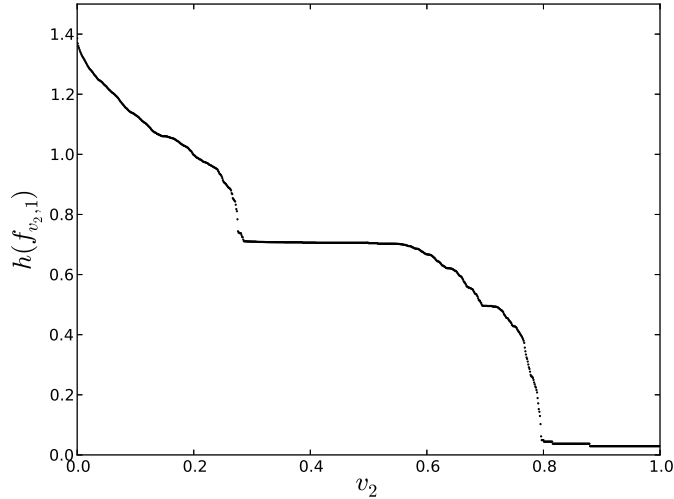}\caption{Plot of $h(f_{v_{2},1})$ vs $v_{2}$, $0\leq v_{2}<1$
($\varepsilon=10^{-4},\Delta v_{2}=0.001$).}%
\label{fig3a}%
\end{figure}

Table 3 displays the performance of the new algorithm as compared to the old
one when computing $h(f_{0.7,1})$. Also this time $t_{orig}>t_{simp}$ for
$\varepsilon=10^{-d}$, $4\leq d\leq7$ (as in Table 1 and 2). Furthermore, we
obtain two exact decimal positions of the topological entropy, $h(f_{0.7,1}%
)=0.48...$, with both algorithms and the $\varepsilon$'s considered.

\begin{table}[h]
\centering
\begin{tabular}
[c]{|c|c|c|c|c|c|c|}\hline
& $h_{orig}$ & $n_{orig}$ & $t_{orig}$ & $h_{simp}$ & $n_{simp}$ & $t_{simp}%
$\\\hline
$\varepsilon=10^{-4}$ & 0.494586 & 135 & 0.304254 & 0.49579 & 147 &
0.251014\\\hline
$\varepsilon=10^{-5}$ & 0.48545 & 426 & 2.920753 & 0.48583 & 464 &
2.465768\\\hline
$\varepsilon=10^{-6}$ & 0.482554 & 1345 & 28.715971 & 0.482675 & 1465 &
24.594864\\\hline
$\varepsilon=10^{-7}$ & 0.481637 & 4250 & 290.400729 & 0.481675 & 4630 &
256.136796\\\hline
\end{tabular}
\caption{Comparison of performances when computing $h(f_{0.7,1})$.}%
\label{table3}%
\end{table}

Finally, Fig. 7 depicts the values of $h(f_{v_{2},v_{3}})$ for $v_{2}+0.3\leq
v_{3}\leq1$, $\varepsilon=10^{-4}$, and $\Delta v_{2},\Delta v_{3}=0.01$.

\begin{figure}[h]
\centering
\includegraphics[width=10cm]
	{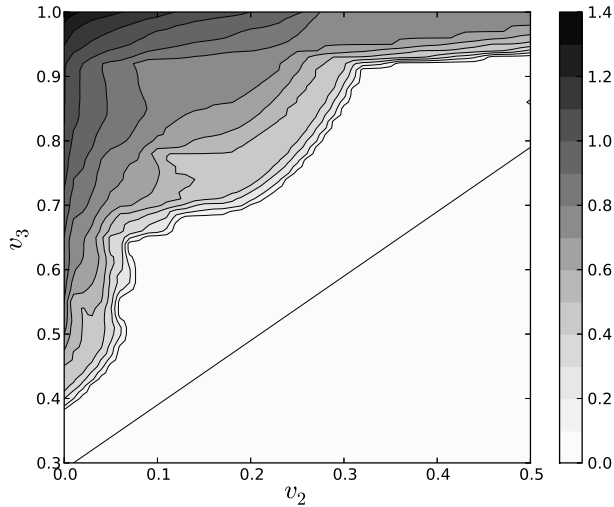}\caption{Level sets of $h(f_{v_{2},v_{3}})$ vs $v_{2},v_{3}$,
$v_{2}+0.3\leq v_{3}\leq1$ ($\varepsilon=10^{-4},\Delta v_{1}=\Delta
v_{2}=0.01$).}%
\label{fig3b}%
\end{figure}

A concluding observation. As anticipated in the remark R2 of Sect.
\ref{section7} and illustrated in the Tables 1-3, the values of $h_{simp}$
converge from above with ever more computation loops (or smaller values of the
parameter $\varepsilon$). This property follows for $h_{orig}$ from \cite[Thm.
4.2.4]{Alseda}.

\section{Conclusion}

We provided in Thm. 1, Eqn. (\ref{Theorem}), an analytical formula to
calculate the topological entropy of a multimodal map $f$. The peculiarity of
Eqn. (\ref{Theorem}), as compared to similar formulas (see (\ref{1.1}%
)-(\ref{1.4})), is that it involves the min-max sequences of $f$. Min-max
sequences generalize kneading sequences in that they contain additional,
geometric information (about the extrema structure of $f^{n}$, $n\geq1$) but
with no computational penalty. We also discussed in Sect.~\ref{section6} the
relationship between the value of $h(f)$ and the divergence rate of
$\log(1-\Sigma_{n})$ (or convergence rate of $\Sigma_{n}$ to $1$). It turns
out that both are related in the way stated in Thms. 2 and 3. A practical
offshoot of Thm. 1 is the algorithm of Sect. \ref{section7} to compute $h(f)$.
This algorithm is a simplified version of a previous one derived in
\cite{Amigo2012}. The performances of both algorithms were compared in Sect. 8
using uni-, bi-, and trimodal maps. In view of the results summarized in
tables 1 to 3, the original algorithm seems to perform better in the unimodal
case, while the opposite occurs in the bimodal and trimodal cases.

\bigskip

Acknowledgements. We thank very much the referees of this paper for their
constructive criticism. We are also grateful to Jos\'{e} S. C\'{a}novas and
Mar\'{\i}a Mu\~{n}oz Guillermo (Universidad Polit\'{e}cnica de Cartagena,
Spain), and V\'{\i}ctor Jim\'{e}nez (Universidad de Murcia, Spain) for
clarifying discussions. This work was financially supported by the Spanish
\textit{Ministerio de Econom\'{\i}a y Competitividad}, grant MTM2012-31698.

\end{document}